\documentstyle[12pt]{article}
\normalbaselines
\begin{document}
\newcommand{\C}{{\bf C}}
\newcommand{\Z}{{\bf Z}}
\newcommand{\R}{{\bf R}}
\newtheorem{th}{Theorem}
\newtheorem{lemma}{Lemma}
\newtheorem{prop}{Proposition}
\newtheorem{add}{Addendum}
\newtheorem{rem}{Remark}

\title{The complex separation and extensions of Rokhlin congruence for curves
on surfaces}
\author{Grigory Mikhalkin}
\date{}
\maketitle

\begin{abstract}
The subject of this paper is the problem of arrangement of real algebraic
curves on real algebraic surfaces. In this paper we extend Rokhlin,
Kharlamov-Gudkov-Krakhnov and Kharlamov-Marin congruences for curves on
surfaces and give some applications of this extension.

 For some pairs consisting of a surface and a curve on this surface
(in particular for M-pairs) we introduce a new structure --- the complex
separation that is separation of the complement of curve into two surfaces.
In accordance with Rokhlin terminology the complex separation is
a complex topological characteristic of real algebraic varieties.
The complex separation is similar to complex orientations introduced by
O.Ya.Viro (to the absolute complex orientation in the case when a curve is
empty and to the relative complex orientation otherwise).
In some cases we calculate the complex separation of a surface
(for example in the case when surface is the double branched covering of
another surface along a curve). With the help of these calculations
applications of the extension of Rokhlin congruence gives some new restrictions
for complex orientations of curves on a hyperboloid.
\end{abstract}

\section{Introduction}
\subsection{Rokhlin and Kharlamov-Gudkov-Krakhnov congruences
for curves on surfaces}
If $A$ is a real curve of odd degree on real projective surface $B$
then $A$ divides $B$ into two parts $B_+$ where polynomial determining $A$
is non-negative and $B_-$ where this polynomial is non-positive.
V.A.Rokhlin \cite{R} proved the congruence for the Euler characteristic
of $B_-$ under such strong hypotheses that they follow that $B$ is an
M-surface,
$A$ is an M-curve and $B_+$ is contained in a connected component of $B$.
V.M.Kharlamov \cite{Kh}, D.A.Gudkov and A.D.Krakhnov \cite{GK} proved a
relevant
congruence that makes sense sometimes even if $A$ is not M- but (M-1)-curve,
but the hypotheses for $B$ and $B_+$ were not weakened.

\subsection{Description of the paper}
The results of the paper make sense in the case when a pair consisting
of a surface and a curve in this surface is of characteristic type
(for definition see section \ref{not}).
For pairs of characteristic type we introduce a complex separation of the
complement of the curve in the surface that is a new complex topological
characteristic of pairs consisting of a surface and a curve in this surface.
If the curve is empty then the complex separation is a new complex topological
characteristic of surfaces relevant to the complex orientation of surfaces
introduced by O.Viro \cite{V1}.
The main theorem is formulated with the help of the complex separation,
this theorem is a generalization of Rokhlin and Kharlamov-Gudkov-Krakhnov
congruences for curves in surfaces.
The main theorem gives nontrivial restrictions even for curves of odd degree
on some surfaces.

The paper contains also some applications of the main theorem.
We prove a congruence modulo 32 for Euler characteristic of real connected
surface of characteristic type.
We prove some new congruences for curves on a hyperboloid.
We give a direct extension of Rokhlin and Kharlamov-Gudkov-Krakhnov
congruences for curves on projective surfaces, there we avoid both of
the hypotheses that $B$ is an M-surface and $B_+$ is contained in a connected
component of $B$.
We apply the main theorem to the classification of curves of low degrees
on an ellipsoid.
In particular, we get a complete classification of flexible curves
on an ellipsoid of bidegree (3,3) (the notion of flexible curve
is analogous to one introduced by O.Viro \cite{V} for plane case).

One can see that the fact that the theorem can be applied to curves of odd
degree follows that this theorem can not be proved in the Rokhlin approach
using the double covering of the complexification of the surface branched along
the curve (since there is no such a covering for curves of odd degree).
We use the Marin approach \cite{Marin}. All results of this paper apply to
flexible curves as well as to algebraic.

The author is indebted to O.Ya.Viro for his attention to the paper and
consultations.

\section{Notations and the statement of the main theorem}
\label{not}

Let $\C B$ be a smooth oriented 4-manifold such that its first $\Z_2$-Betti
number is zero equipped with an involution $conj$ such that the set $\R B$ of
its fixed points is a surface.
Let $\C A$ be a smooth surface in $\C B$ invariant under $conj$ and such that
the intersection of $\C A$ and $\R B$ is a curve.
These notations are inspired by algebraic geometry.

It is said that $A$ is of even degree if $\C A$ is $\Z_2$-homologous to zero
in  $\C B$ and that $A$ is of odd degree otherwise.
It is said that curve $A$ is of type I if $\R A$ is $\Z_2$-homologous to zero
in $\C A$ and that $A$ is of type II otherwise.
It is said that surface $B$ is of type I$abs$ if $\R B$ is $\Z_2$-homologous
to zero in $\C B$
If $B$ is a real projective plane then it is said that $B$ is of type I$rel$
if $\R B$ is $\Z_2$-homologous to a plane section of $\C B$.
We shall say that pair $(B,a)$ is of characteristic type if the sum of
$\R B$ and $\C A$ is Poincar\'{e} dual to the second Stiefel-Whitney class
of $\C B$.
We shall say that surface $B$ is of characteristic type if $\R B$ is a
characteristic surface in $\C B$.

Let $b_*$ denote the total $\Z_2$-Betti number.
It is said that manifold $\C X$ equipped with involution $conj$ is an
(M-$j$)-manifold if $b_*(\R X)+2j=b_*(\C X)$ where $\R X$ is the fixed point
set
of $conj$.
One can easily see that Smith theory follows that $j$ is a nonnegative integer
number.
Let $\sigma(M)$ denote the signature of oriented manifold $M$.;
$D_M:H^*(M;\Z_2)\rightarrow H_*(M;\Z_2)$ denote Poincar\'{e} duality operator;
$[N]\in H_*(M;\Z_2)$ denote $\Z_2$-homology class of submanifold $N\subset M$.
Let $e_A=[\C A\circ\C A]_{\C B}$ denote normal Euler number of $\R B$ in $\C
B$.
If $B_\epsilon$ is a surface contained in $\R B$ and such that
$\partial B_\epsilon=\R A$ that we shall denote by $e_{B_\epsilon}$
the obstruction to extending of line bundle over $\R A$ and normal to $\R A$
in $\C A$ to the line bundle over $B_\epsilon$ normal in $\C B$ to $\R B$
evaluated on the twisted fundamental class $[B_\epsilon,\partial B_\epsilon]$
and divided by 2.
One can see that if $(B,A)$ is a nonsingular pair consisting of
an algebraic surface and algebraic curve then
$e_{B_\epsilon}=-\chi(B_\epsilon)$.
Let $\beta(q)$ denote the Brown invariant of $\Z_4$-valued quadratic form $q$.

\begin{th}
\label{main}
If $(B,A)$ is of characteristic type then there is a natural separation of
$\R B-\R A$ into surfaces $B_1$ and $B_2$ such that
$\partial B_1=\partial B_2=\R A$
defined by the condition that $\C A/conj\cup B_j$ is a characteristic surface
in $\C B/conj (j=1,2)$.
There is a congruence for the Guillou-Marin form $q_j$ on
$H_1(\C A/conj\cup B_j;\Z_2)$
\begin{displaymath}
e_{B_j}\equiv\frac{e_{\R B}+\sigma(\C B)}{4}-\frac{e_A}{4}-\beta(q_j)\pmod{8}
\end{displaymath}
\end{th}
\begin{add}
Let $q_j|_{H_1(\R A;\Z_2)}=0$
\begin{itemize}
\begin{description}
\item[a)] If $A$ is an M-curve then
$e_{B_j}\equiv\frac{e_{\R B}+\sigma(\C B)}{4}-\frac{e_A}{4}-\beta_j\pmod{8}$
\item[b)] If $A$ is an (M-1)-curve then
$e_{B_j}\equiv\frac{e_{\R B}+\sigma(\C
B)}{4}-\frac{e_A}{4}-\beta_j\pm1\pmod{8}$
\item[c)] If $A$ is an (M-2)-curve and
$e_{B_j}\equiv\frac{e_{\R B}+\sigma(\C B)}{4}-\frac{e_A}{4}-\beta_j+4\pmod{8}$
then $A$ is of type I
\item[d)] If $A$ is of type I then
$e_{B_j}\equiv\frac{e_{\R B}+\sigma(\C B)}{4}-\frac{e_A}{4}\pmod{4}$
\end{description}
\end{itemize}
where $\beta_j$ is  the Brown invariant of the restriction
$q_j|_{H_1(B_j;\Z_2)}$
\end{add}

\subsection{Remark}
Some of components of $\R A$ can be disorienting loops in $\C A$
(it is easy to see that number of such components is even).
If $\alpha$ is some 1-dimensional $\Z_2$-cycle in $B_j$ that is a boundary
of some 2-chain $\beta$ in $\R B$ containing an even number of
disorienting $\C A$ components of $\R A$ then $q_j(\alpha)=0$; if such a number
is odd
then $q_j(\alpha)=2$.
It follows that if $(B,A)$ is of characteristic type then the number
of disorienting $\C A$ components of $\R A$ in each component of $\R B$
is even.

\section{The proof of Theorem 1 and Addendum 1}
\subsection{Calculation of the characteristic class of $\C B/conj$}
It is not difficult to see the formula for the characteristic classes of
double branched covering:
if $\pi:Y\rightarrow X$ is a double covering branched along $Z$ then
$$w_2(Y)=\pi^*w_2(X)+D^{-1}_Y[Z]$$
Applying this formula to $p:\C B\rightarrow \C B/conj$ we get
$$tr(Dw_2(\C B/conj))=Dw_2(\C B)+[\R B]$$
where $tr:H_2(\C B/conj;\Z_2)\rightarrow H_2(\C B;\Z_2)$ is transfer
(i.e. the inverse Hopf homomorphism to $P$).
It is easy to see that transfer can be decomposed as the composition
$$H_2(\C B/conj;\Z_2)\stackrel{k}{\rightarrow}H_2(\C B/conj,\R B;\Z_2)
\stackrel{h}{\rightarrow}H_2(\C B;\Z_2)$$
where $k$ is an inclusion homomorphism and $h\circ k=tr$.
To prove that $h$ is a monomorphism we use the Smith exact sequence
(see e.g. \cite{W}):
$$H_3(\C B/conj,\R B;\Z_2)\stackrel{\gamma_3}{\rightarrow}
H_2(\R B;\Z_2)\oplus H_2(\C B/conj,\R B;\Z_2)\stackrel{\alpha_2}{\rightarrow}
H_2(\C B;\Z_2)$$
In this sequence the first component of $\gamma_3$ is equal to
the boundary homomorphism $\partial$ of pair $(\C B/conj,\R B)$;
$\partial$ is a monomorphism since $H_3(\C B/conj;\Z_2)=0$
(since $\C B$ and therefore $\C B/conj$ are simply connected).
It means that no element of type $(0,x)\in H_2(\R B;\Z_2)\oplus
H_2(\C B/conj,\R B;\Z_2), x\neq 0$ is contained in $Im\gamma_3$
and therefore the restriction of $\alpha_2$ to $H_2(\C B/conj,\R B;\Z_2)$
is a monomorphism and thus $H$ is a monomorphism.
Now if $[\C A]=Dw_2(\C B)+[\R B]$ then
$$Dw_2(\C B/conj)=[\C A/conj]\in H_2(\C B/conj,\R B;\Z_2)$$
The exactness of homology sequence of pair $(\C B/conj,\R B)$ follows
since $H$ is a monomorphism that there exists a surface $B_1\subset\R B$
such that $\partial B_1=\R A$ and $W_1=B_1\cup\C A/conj$ is dual to
$w_2(\C B/conj)$. Let $B_2$ be equal to $Cl(\R B-B_1)$
Surface $\R B$ is $\Z_2$-homologous to zero in $\C B/conj$ since $\C B$ is
a double covering of $\C B/conj$ branched along $\R B$.
It follows that $W_2=B_2\oplus\C A/conj$ is also a surface dual to
$w_2(\C B/conj)$.

Note that the separation of $\R B$ into $B_1$ and $B_2$ is unique since
$$dim(ker(in_*:H_2(\R B;\Z_2)\rightarrow H_2(\C B/conj)))=1$$
as it follows from exactness of the Smith sequence.
Indeed, since $H_3(\C B/conj;\Z_2)=0$ then this dimension is equal to the
dimension of $H_3(\C B/conj,\R B;\Z_2)$;
the dimension of $H_3(\C B/conj,\R B;\Z_2)$ is equal to 1 since
$\gamma_4:H_4(\C B/conj,\R B;\Z_2)\rightarrow {0}\oplus H_3(\C B/conj,\R B;
\Z_2)$ is an isomorphism.

\subsection{Proof of the congruence for $\chi(B_j)$}
\label{prth}
Note that since $W_j,j\in \{1,2\}$ is a characteristic surface in $\C B/conj$
and $\C B$ is simply connected
we can apply the Guillou-Marin congruence to pair $(\C B/conj,W_j)$
$$\sigma(\C B/conj)\equiv W_j\circ W_j+2\beta (q_j)\pmod{16}$$
where $q_j:H_1(W_j;\Z_2)\rightarrow\Z_4$ is the quadratic form associated
to the embedding of $W_j$ into $\C B/conj$ (see \cite{GM}).
Similar to the calculations in \cite{Marin} we get that
$$W_j\circ W_j=\frac{e_A}{2}-2\chi(B_j)$$
The Atiyah-Singer-Hirzebruch formula follows that
$$\sigma(\C B/conj)=\frac{\sigma(\C B)-\chi(\R B)}{2}$$
Combining all this we get
$$\chi(B_j)\equiv\frac{e_A}{4}+\frac{\chi(\R B)-\sigma(\C B)}{4}+
\beta(q_j)\pmod{8}$$
Now if $q|_{H_1(\R A;\Z_2)}=0$ then additivity of the Brown invariant
(see \cite{KV}) follows that
$$\beta(q_j)=\beta(q|_{H_1(\C A/conj;\Z_2)})+\beta_j$$
It is easy to see that if $A$ is an (M-$j$)-curve then
$rkH_1(\C A/conj;\Z_2)=j$.
Points a) and b) of the addendum immediately follow from this.
To deduce points c) and d) of the addendum note that
$q|_{H_1(\C A/conj;\Z_2)}$ is even iff
$\C A/conj$ is an orientable surface iff $A$ is of type I (cf. \cite{Marin},
\cite{KV}).

\section{Some applications of Theorem 1}
\subsection{The case $\C A$ is empty; congruences for surfaces}
\label{pov}
\subsubsection{}
\label{surface}
If $Dw_2(\C B)=[\R B]$ then there is defined a complex separation of $\R B$
into two closed surfaces $B_1$ and $B_2$;
there is defined a $\Z_4$-quadratic form $q$ on $H_1(\R B;\Z_2)=
H_1(B_1;\Z_2)\oplus H_1(B_2;\Z_2)$ equal to sum of Guillou-Marin forms
of $B_1$ and $B_2$ which are characteristic surfaces in $\C B/conj$ and
\begin{displaymath}
\chi(B_j)\equiv\frac{\chi(\R B)-\sigma(\C B)}{4}+\beta(q|_{H_1(B_j;\Z_2)})
\pmod{8}
\end{displaymath}

\subsubsection{}
\label{pusto}
If $Dw_2(\C B)=[\R B]$ and $B_j$ is empty for some $j$
(that is evidently true if $\R B$ is connected) then
$$\chi(\R B)\equiv\sigma(\C B)\pmod{32}$$

\subsubsection{}
If $Dw_2(\C B)=[\R B]$ then $$\chi(\R B)\equiv\sigma(\C B)\pmod{8}$$
{\em\underline{Proof}} It follows from an easy observation that $\chi(B_j)
\equiv\beta(q|_{H_1(B_j;\Z_2)})\pmod{2}$

\subsubsection{Remark}
According to O.Viro \cite{V} in some cases one can define some more complex
topological characteristics on $\R B$.
Namely, If $B$ is of type I$abs$ then $\R B$ possesses two special
reciprocal orientations (so-called semi-orientation) and a special spin
structure.
If $Dw_2(\C B)=[\R B]$ then $\R B$ possesses the special $Pin_-$-structure
corresponding to Guillou-Marin form $q_{\C B}:H_1(\R B;\Z_2)\rightarrow\Z_4$
of surface $\R B$ in $\C B$.
The complex separation is a new topological characteristic
for surfaces and \ref{surface} may be interpreted as a formula for this
characteristic.
Quadratic form $q$ is not a new complex topological characteristic.

\subsubsection{}
Form $q$ is equal to $q_{\C B}$ in the case when these forms are defined (i.e.
when
$Dw_2(\C B)=[\R B]$)

To prove this one can note that the index of a generic membrane in $\C B$
bounded by curves in $\R B$ differs from the index of the image of this
membrane in $\C B/conj$ by number of intersection points of this membrane
and $\R B$.

\subsubsection{Remark}
If $B$ is a complete intersection in the projective space of hypersurfaces
of degrees $m_j,j=1,\ldots,s$ then the condition that $Dw_2(\C B)=[\R B]$
is equivalent to the condition that $B$ is of type I$abs$ in the case
when $\sum^{s}_{j=1}m_j\equiv 0\pmod{2}$ and to the condition that $B$ is
of type I$rel$ in the case when $\sum^{s}_{j=1}m_j\equiv 1\pmod{2}$.

\subsection{Calculations for double coverings}
\label{vych}
We see that to apply Theorem 1 one needs to be able to calculate the complex
separation and the corresponding Guillou-Marin form.
We calculate them in some cases in this subsection.

By the semiorientation of manifold $M$ we mean a pair of reciprocal
orientations of $M$ (note that this notion is nontrivial only for non-connected
manifolds).
It is easy to see that two semiorientations determine a separation of $M$;
this separation is a difference of two semiorientations,
namely, two components of $M$ are of the same class of separation iff the
restrictions of the semiorientations on these components are the same.

Let $\C B$ be the double covering of surface $\C X$
branched along curve $\C D$
invariant under the complex conjugation $conj_X$ in $\C X$.
Suppose $\C D$ is of type I.
Let $\C D_+$ be one of two components of $\C D-\R D$.
Then $\R D$ possesses a special semiorientation called the complex
semiorientation (see \cite{R1}).
The invariance of $\C D$ under $conj_X$ follows that $conj_X$ can be lifted
in two different ways into an involution of $\C B$.
Let $conj_B$ be one of these two lifts.
Let $X_-=p(\R B))$, $X_+=\R X-int(X_-)$ where $p$ is the covering map.

\subsubsection{}
\label{calor}
Suppose that $Dw_2(\C B)=[\R B]$.
Then for every component $C$ of $X_+$ the complex separation of $\partial C$
induced from complex separation of $\R B$ via $p$ (namely, two circles of
$\partial C$ are of the same class of separation iff the lie in the image
under $p$ of the same class of the separation of $\R B$)
is equal to the difference of the semiorientations on $\partial C$
induced by the complex semiorientation of $\R D$ and the unique (since $C$ is
connected) semiorientation of $C$.
In particular $C$ is orientable.

{\em\underline{Proof}} Let $\alpha$ and $\beta$ be components of $\partial C$.
Consider two only possible cases: the first case when the semiorientations
induced from $\R D$ and $C$ are equal on $\alpha$ and $\beta$ and the second
case when they are different (see fig.1, arrows indicate one of two
orientations induced by the complex semiorientation of $\R D$).
Choose a point $Q_{\alpha}$ in $\alpha$ and $Q_{\beta}$ in $\beta$.
Connect these points by a path $\gamma$ inside $C$ and by a path $\delta$
inside $\C D_+$ (not visible on the picture) without self-intersection points.
It is easy to see that there exists a disk $F'\subset\C X$ bounding loop
$\gamma\delta$ and such that the interior of $F'$ does not intersect $\C D$.
Set $F$ to be equal to $F'\cup conj_X F'$.
Then $p^{-1}(F)$ gives an element of $H_2(\C B)$ (the construction of this
element was suggested by O.Viro \cite{V3}).
Since $p^{-1}(F)$ is invariant under $conj_B$ it gives an element in
$H_2(\C B/conj_B;\Z_2)$, say $f\in H_2(\C B/conj_B;\Z_2)$.
Let us calculate the self-intersection number of $f$.
It is easy to see that because of symmetry the self-intersection number of $f$
in $\C B/conj_B$ is equal to the self-intersection number of $\gamma\delta$
in $C\cup\C D_+$.
The definition of complex semiorientation follows that the self-intersection
number of $\gamma\delta$ in $C\cup\C D_+$ (and therefore the self-intersection
number of $f$) is equal to zero in the first case and to one in the second
case.

\subsubsection{Remark}
In the case when $\R X$ is connected \ref{calor}
completely determines the complex separation of $\R B$.

\subsubsection{(O.Viro [11])}
\label{vychv}
If $\lambda$ is a loop in $\R B$ such that $p(\lambda)=\partial(G)$,
where $G\subset\R X$ then
$$q_{\C B}(\lambda)\equiv 2\chi(G\cap X_+)\pmod{4}$$

\subsubsection{(O.Viro [11])}
\label{vychq}
Suppose $\gamma$ is a path in $X_-$ connecting points $Q_{\alpha}$ and
$Q_{\beta}$ of components $\alpha$ and $\beta$ of $\R D$ respectively.
Then $q_{\C B}(p^{-1}(\gamma))=0$ if the intersection numbers of $\gamma$ with
$\alpha$ and
$\beta$ are of opposite sign (case 1 of fig.1) and
$q_{\C B}(p^{-1}(\gamma))=2$ otherwise (case 2 of fig.1)

\subsection{New congruences for complex orientations of curves on a
hyperboloid}
In this subsection we apply the results of \ref{pov} and \ref{vych} to
double branched coverings over simplest surfaces of characteristic type,
a plane and a hyperboloid.
To state congruences it is convenient to use the language of integral calculus
based on Euler characteristic developed by O.Viro \cite{V2}.
Let $\R A$ be a curve of type I in the connected surface $\R X$.
We equip $\R A$ with one of two complex orientations and fix $X_{\infty}$
-- one of the components of $\R X-\R A$.
If $\R X-X_{\infty}$ is orientable and $\R A$ is $R$-homologous to zero
for some ring $R$ of coefficient coefficient
then there is defined function $ind_R:\R X-\R A\rightarrow R$ equal to zero
on $X_{\infty}$ and equal to the $R$-linking number with oriented curve $\R A$
in $\R X-X_{\infty}$ otherwise.
It is easy to see that $ind_R$ is measurable and defined almost everywhere on
$\R X$ with respect to Euler characteristic.
Evidently, the function ${ind_R}^2:\R X-\R A\rightarrow R\otimes R$ does
not depend on
the ambiguity in the choice of one of two complex orientation of $\R A$.

\subsubsection{}
Consider the case when $X=P^2$, $R=\Z$.
Let $A$ be a plane nonsingular real curve of type I given by polynomial $f_A$
of degree $m=2k$.
Then $\R A$ is $\Z_2$-homologous to zero.
Let $X_{\infty}$ be the only nonorientable component of $\R P^2-\R A$.
Without loss of generality suppose that $f_A|_{X_{\infty}}<0$.
Define $X_{\pm}$ to be equal to $\{y\in\R P^2|\pm f_A(y)\ge 0\}$.
Let $p:\C B\rightarrow\C P^2$ be the double covering of $\C P^2$ branched
along $\C A$ (note that such a covering exists and is unique since $m$ is even
and $\C P^2$ is simply connected).
Let $conj_B:\C B\rightarrow\C B$ be the lift of $conj:\C P^2\rightarrow\C P^2$
such that $p(\R B)=X_-$ (as it is usual we denote $Fix(conj_B$ by $\R B$).
Lemmae 6.6 and 6.7 of \cite{W} immediately follow that $Dw_2(\C B)=\R B$.
Therefore we can apply \ref{surface} and \ref{calor}.
Proposition \ref{calor} follows that $p^{-1}(cl({ind_{\Z_2}}^{-1}(2+4\Z)))$ is
equal to one of two surfaces of the complex separation of $\R B$.
Set $B_1$ to be equal to $p^{-1}(cl({ind_{\Z}}^{-1}_{\Z_2}(2+4\Z)))$,
note that $B$ is orientable since $ind_{\Z}|_{X_{\infty}}=0$.

\subsubsection{Lemma}
$$\beta(q|_{H_1(B_1;\Z_2)})\equiv
4\chi({ind_{\Z}}^{-1}(3+8\Z)\cup{ind_{\Z}}^{-1}(-3+8\Z))
\pmod{8}$$
{\em\underline{Proof}}The boundary of each component $C$ of $p(B_1)$ has one
exterior oval and some interior ovals (with respect to $C$).
We call an interior oval of $C$ $C$-positive if its complex orientation and
the complex orientation of the exterior oval can be extended to some
orientation
of $C$ and we call it $C$-negative otherwise.
Proposition \ref{vychq} follows that $\beta(q|_{H_1(B_1;\Z_2)})$ is equal to
twice the sum of values of $q$ on on all $C$-negative ovals modulo 8.
Note that $ind_{\Z}^{-1}(\pm 3+8\Z)$ is just the part of $X_+$ lying inside
odd number of $C$-negative ovals.
The lemma follows now from \ref{vychv}.

\subsubsection{}
The application of \ref{surface} gives that
$$\chi(B_1)\equiv\frac{\chi(\R B)-\sigma(\C B)}{4}+\beta(q|_{H_1(B_1;\Z_2)})
\pmod{8}$$
Note that $p(\R B)={ind_{\Z}}^{-1}(2\Z)$, $\chi({ind_{\Z}}^{-1}(2\Z))=
1-\chi(1+2\Z)$
and $\sigma(\C B)=2-2k^2$.
We get
$$4\chi({ind_{\Z}}^{-1}(2+4\Z))\equiv 1-\chi({ind_{\Z}}^{-1}(1+2\Z))-1+k^2+
8\chi({ind_{\Z}}^{-1}(\pm 3+8\Z))\pmod{16}$$
We can reformulate this in integral calculus language
$$\int_{\R P^2}{ind_{\Z}}^2d\chi\equiv k^2\pmod{16}$$
Thus for projective plane we get nothing new but the reduction modulo 16
of the Rokhlin congruence for complex orientation \cite{R1}
$$\int_{\R P^2}{ind_{\Z}}^2d\chi=k^2$$

\subsubsection{}
Consider now the case when $X=P^1\times P^1$.
Let $A$ be a nonsingular real curve of type I in $P^1\times P^1$ of bidegree
$(d,r)$, i.e. the bihomogeneous polynomial $f_A$ determining $A$ is of bidegree
$(d,r)$ where $d$ and $r$ are even numbers.
Let $X_{\infty}$ be a component of $\R P^2-\R A$, $X_{\pm}=
\{y\in\R P^1\times\R P^1|\pm f_A(y)\ge 0\}$.
Suppose without loss of generality that $f_A|_{X_{\infty}}<0$.
Let $p:\C B\rightarrow\C P^2$ be the double covering of $\C P^2$ branched
along $\C A$ and let $conj_B:\C B\rightarrow\C B$ be the lift of
$conj:\C P^1\times\C P^1\rightarrow\C P^1\times\C P^1$ such that
$p(\R B)=X_-$.
Nonsingularity of $A$ follows that all components of $\R A$ non-homologous to
zero are homologous to each other.
Let $e_1,e_2$ form the standard basis of $H_1(\R P^1\times\R P^1)$ and
let $s,t$ be the coordinates in this basis of a non-homologous to zero
component of $\R A$ equipped with such an orientation that $s,t\ge 0$.
If all the components of $\R A$ are homologous to zero then set $s=t=0$.
Then $\R A$ equipped with the complex orientation produces $l'(se_1+te_2)$
in $H_1(\R P^1\times\R P^1)$.
Note that $s$ and $t$ are relatively prime and $l'$ is even since both $d$ and
$r$ are even.

\subsubsection{Lemma}
If $l'\equiv 0\pmod{4}$ then $Dw_2(\C B)=[\R B]$

The proof follows from Lemma 3.1 of \cite{Ma}.

\subsubsection{Lemma}
If $l'\equiv 0\pmod{4}$ then $ind_{\Z_4}$ is defined
and if $sd+tr\equiv 0\pmod{4}$ then
$$\int_{\R P^1\times\R P^1}{ind_{\Z_4}}^2d\chi\equiv \frac{dr}{2}\equiv
0\pmod{8}$$
{\em\underline{Proof}} Note that the condition that $sd+tr\equiv 0\pmod{4}$
is just equivalent to the orientability of $\R B$.
Therefore $\beta(q|_{H_1(B_1;\Z_2)})\equiv 0\pmod{4}$.
Further arguments are similar to the plane case,
we skip them.

\subsubsection{Remark}
The traditional way of proving of formulae of complex orientations for
curve on surfaces (see \cite{Z}) adjusted to the case when the real curve
is only $\Z_4$-homologous to zero gives only congruence
$$\int_{\R P^1\times\R P^1}{ind_{\Z_4}}^2d\chi\equiv \frac{dr}{2}\pmod{4}$$
The \ref{l=4} shows the worthiness of modulo 4 in this congruence.
If $l'\equiv 0\pmod{8}$ then the traditional way gives that
$$\int_{\R P^1\times\R P^1}{ind_{\Z_8}}^2d\chi\equiv \frac{dr}{2}\pmod{8}$$

\subsubsection{}
\label{l=4}
If $l'\equiv 4\pmod{8}$ and $sd+tr\equiv 2\pmod{4}$
then
$$\int_{\R P^1\times\R P^1}{ind_{\Z_4}}^2d\chi\equiv \frac{dr}{2}+4\pmod{8}$$
{\em\underline{Proof}} Form $q|_{H_1(B_1;\Z_2)}$ is cobordant to the sum of
a form on an orientable surface and some forms on Klein bottles.
The condition that $l'\equiv 4\pmod{8}$ is equivalent to the condition
that the number of forms on Klein bottles non-cobordant to zero  is odd.
Thus
$$\beta(q|_{H_1(B_1;\Z_2)})\equiv 2\pmod{4}$$

\subsubsection{}
If $l'\equiv 0\pmod{8}$ and $sd+tr\equiv 0\pmod{4}$ then
$$\int_{\R P^1\times\R P^1}{ind}_{\Z_8}^{2}d\chi\equiv \frac{dr}{2}\pmod{16}$$
The proof is similar to the plane case.

\subsubsection{Addendum, new congruences for the Euler characteristic of $B_+$
for curves on a hyperboloid}
\label{b10}
Let $d\equiv r\equiv 0\pmod{2}$ and $\frac{d}{2}t+\frac{r}{2}s+s+t\equiv
1\pmod{2}$.
\begin{itemize}
\begin{description}
\item[a)] If $A$ is an M-curve then $\chi(B_+)\equiv\frac{dr}{2}\pmod{8}$
\item[b)] If $A$ is an (M-1)-curve then $\chi(B_+)\equiv\frac{dr}{2}\pm
1\pmod{8}$
\item[c)] If $A$ is an (M-2)-curve and $\chi(B_+)\equiv\frac{dr}{2}+4\pmod{8}$
then $A$ is of type I
\item[d)] If $A$ is of type I then $\chi(B_+)\equiv 0\pmod{4}$
\end{description}
\end{itemize}

This theorem follows from Theorem 1 and gives some new restriction on the
topology of the arrangement of real nonsingular algebraic curve of even
bidegree on a hyperboloid with non-contractible branches.
Points a) and b) of \ref{b10} in the case when $\frac{d}{2}t+\frac{r}{2}s\equiv
0\pmod{2}, s+t\equiv 1\pmod{2}$ were proved by S.Matsuoka \cite{Ma1} in another
way (using 2-sheeted branched coverings of hyperboloid).
Point d) of \ref{b10} is a corollary of the modification of Rokhlin formula
of complex orientations for modulo 4 case.

\subsection{The case when $A$ is a curve of even degree on projective surface
$B$}
In this subsection we deduce Rokhlin and Kharlamov-Gudkov-Krakhnov congruences
from Theorem 1 and give a direct generalization of these congruences.

Let $B$ be the surface in $P^q$ given by the system of equations
$P_j(x_0,\ldots,x_q)=0,j=1,\ldots,s-1$,
let $A$ be the (M-$k$)-curve given by the system of equations
$P_j(x_0,\ldots,x_q)=0,j=1,\ldots,s-1$,
where $P_j$ are homogeneous polynomials with real coefficients,
$deg P_j=m_j$, $s=q-1$.
Suppose $(B,A)$ is a non-singular pair, $\R A\neq\emptyset$ and $m_s$ is even.
Denote $B_+=\{x\in\R B|P_s(x)\geq 0\},B_-=\{x\in\R B|P_s(x)\leq 0\}$,
$$d=rk(in^{B_+}_*:H_1(B_+;\Z_2)\rightarrow H_1(\R B;\Z_2)),
e=rk(in^{\R A}_*:H_1(\R A;\Z_2)\rightarrow H_1(\R B;\Z_2))$$
Set $c$ to be equal to the number of non-contractible in $\R P^q$ components
of $\R B$ not intersecting $\R A$.

Let us reformulate Rokhlin and Kharlamov-Gudkov-Krakhnov congruences for
curve on surfaces in the form convenient for generalization and correcting
the error in \cite{R}.

\subsubsection{(Rokhlin [1], Kharlamov [2], Gudkov-Krakhnov [3])}
\label{RKGK}
Suppose $B$ is an $M$-surface, $e=0$, $B_+$ is contained in one component
of $\R B$ and in the case when $m_s\equiv 0\pmod{4}$ suppose in addition that
$c=0$.
\begin{itemize}
\begin{description}
\item[a)] If $d+k=0$ then
$$\chi(B_+)\equiv\frac{m_1\ldots m_{s-1}m_s^2}{4}\pmod{8}$$
\item[b)] If $d+k=1$ then
$$\chi(B_+)\equiv\frac{m_1\ldots m_{s-1}m_s^2}{4}\pm 1\pmod{8}$$
\end{description}
\end{itemize}

Indeed, it is easy to see that the hypothesis of \ref{RKGK} without
the condition on $c$ is equivalent to the hypotheses of corresponding theorems
in \cite{R}, \cite{Kh} and \cite{GK},
to reformulate them it is enough to apply the Rokhlin congruence for
M-surfaces.

\subsubsection{Remark (on an error in [1])}
Point 2.3 of \cite{R} contains a miscalculation of characteristic class $x$
of the restriction to $B_-$ of the double covering of $\C B$ branched along
$\C A$.
In \cite{R} it is claimed that if $m_s\equiv 2\pmod{4}$ then $x=w_1(B_{-}-A)$.
It led to the omission of the condition on $c$ in both Rokhlin and
Kharlamov-Gudkov-Krakhnov congruences.

\subsubsection{Correction of the error in [1]}
If $m_s\equiv 2\pmod{4}$ then $x=in^*\alpha$ where $in$ is the inclusion of
$B_-$
into $\R P^q$ and $\alpha$ is the only non-zero element of $H^1(\R P^q;\Z_2)$

{\em\underline{Proof}} Consider $E$ -- the auxiliary surface in $P^q$
given by equation $P_s(x_0,\ldots,x_q)=0$.
Then the construction of the double covering of $\C P^q$ branched along $\C E$
in the weighted-homogeneous projective space by equation $\lambda^2=
P_s(x_0,\ldots,x_q)$ follows that the characteristic class of the restriction
of the covering to $\R P^q-\R E$ is equal to the restriction of $\alpha$
to $\R P^q-\R E$.
Therefore, the characteristic class of the restriction of the covering
to $B_-$ is equal to $in^*\alpha$.

Besides, it is claimed in \cite{R} that in the case of curves on surfaces
(i.e. $n=1$ in notations of \cite{R}) the endomorphism
$\omega:H_*(B_-,A;\Z_2)\rightarrow H_*(B_-,A;\Z_2)$ of cap-product with
characteristic class $x$ is trivial.
This is true only if either $m_s\equiv 0\pmod{4}$ or each non-contractible
component of $B_-$ contains at least one component of $\R A$.
The condition on $C$ in \ref{RKGK} allows to correct proofs of congruences
of \cite{R}, \cite{Kh} and \cite{GK}.
The author though does not know counter-examples to \ref{RKGK} without
condition on $c$ (without the condition on $c$ point a) of \ref{RKGK} is
equivalent to 3.4 of \cite{R}).

\subsubsection{Direct generalization of Rokhlin and Kharlamov-Gudkov-Krakhnov
congruences}
\label{gen}
Suppose that $B$ is of type I$abs$ in the case when $\sum_{j=1}^{s-1}m_j\equiv
0
\pmod{2}$ and that $B$ is of type I$rel$ in the case when
$\sum_{j=1}^{s-1}m_j\equiv 1\pmod{2}$.
Suppose that $m_s$ is even, $e=0$, $B_+$ is contained in one surface of
the complex separation of $\R B$ and in the case when $m_s\equiv 2\pmod{4}$
suppose in addition that $c=0$.
\begin{itemize}
\begin{description}
\item[a)] If $d+k=0$ then
$$\chi(B_+)\equiv\frac{m_1\ldots m_{s-1}m_s^2}{4}\pmod{8}$$
\item[b)] If $d+k=1$ then
$$\chi(B_+)\equiv\frac{m_1\ldots m_{s-1}m_s^2}{4}\pm 1\pmod{8}$$
\item[c)] If $d+k=2$ and
$$\chi(B_+)\equiv\frac{m_1\ldots m_{s-1}m_s^2}{4}+4\pmod{8}$$
then $A$ is of type I
\item[d)] If $A$ is of type I then
$$\chi(B_+)\equiv\frac{m_1\ldots m_{s-1}m_s^2}{4}\pmod{4}$$
\end{description}
\end{itemize}

\subsubsection{Lemma}
\label{V+=0}
The surface $V_+=\C A/conj\cup B_+$ is $\Z_2$-homologous to zero in
$\C B/conj$

{\em\underline{Proof}} Consider the diagram
\begin{picture}(200,-100)(-50,0)
\put(-10,0){$\C Y$}\put(60,0){${\bf C}B$}
\put(-20,-40){$\C Y/conj_Y$}\put(60,-40){${\bf C}B/conj .$}
\put(10,4){\vector(1,0){46}}
\put(3,-3){\vector(0,-1){25}}
\put(67,-3){\vector(0,-1){25}}
\put(33,8){$p$}
\end{picture}\\[42pt]
where $p:\C Y\rightarrow\C B$ is the double covering over $\C B$ branched
along $\C A$ and $conj_Y$ is such a lift of $conj:\C B\rightarrow\C B$
that $\R Y=\{y\in\C Y|conj_Yy=y\}$ is equal to $p^{-1}(B_-)$.
It is easy to see that this diagram can be expanded to a commutative one
by adding $\phi:\C Y/conj_Y\rightarrow\C B/conj$ where $\phi$ is the double
covering map branched along $V_+$.
It follows that $[V_+]=0\in H_2(\C B/conj)$.

\subsubsection{Remark}
The construction of \ref{V+=0} allows to define the separation of $\R B$
into $B_+$ and $B_-$ for any (not necessarily algebraic) curve
of even degree and invariant with respect to $conj$
in any (not necessarily projective) complex surface $\C B$ equipped with
almost antiholomorphic involution $conj:\C B\rightarrow\C B$ and such that
$H_1(\C B;\Z_2)=0$.
In this case there is the unique double covering over $\C B$ branched along
$\C A$ and two possible lifts of $conj:\C B\rightarrow\C B$.
The images of the fix point sets of these lifts form the desired
separation.

Note that this separation is different from the complex separation.
Let us give an internal definition of this separation.
Let $x,y$ be points in $\R B-\R A$.
Connect these points by path $\gamma$ inside $\C B-\C A$.
If the linking number of the loop $\gamma conj\gamma$ and $\C A$ is equal
to zero then $x$ and $y$ are points of the same surface of the complex
separation, otherwise $x$ and $y$ are points of two different surfaces
of the complex separation.

This remark allows to extend \ref{gen} to the case of flexible curves.

\subsubsection{The proof of 4.4.4}
Note that $Dw_2(\C B)+[\R B]=0$ since $Dw_2(\C B)\equiv [\infty]
\sum_{j=1}^{s-1}m_j$ where $[\infty]\in H_2(\C B;\Z_2)$ is the class of
hyperplane section of $\C B$ and note that $[\C A]=0$ since $m_s$ is even.

We apply \ref{surface}, let $W$ be the surface of the complex separation
not intersecting $B_+$.
Then \ref{V+=0} follows that $W\cup B_+$ is the surface of the complex
separation of $(B,A)$ since if $W$ is a characteristic surface in
$\C B/conj$ then so is $W_+=W\cup V_+$.

Let us prove now that $q_{W_+}|_{H_1(W;\Z_2)}$ is equal to $q_W$ where
$q_{W_+}$ and $q_W$ are Guillou-Marin forms of $W_+$ and $W$ in
$\C B/conj$.
Note that $(q_{W_+}-q_W)(x), x\in H_1(W;\Z_2)$ is equal to the linking number
of $x$ and $V_+$ in $\C B/conj$ that is equal to the linking number of
$x$ and $\C A$ in $\C B$.
It is not difficult to see that the condition on $c$ follows that
the linking number of $x$ and $\C A$ in $\C B$ is equal to zero.

Theorem 1 follows that
$$\chi(B_+)+\chi(W)\equiv\frac{e_A}{4}+\frac{\chi(\R B)-\sigma(\C B)}{4}+
\beta(q_W)+\beta(q_{W_+}|_{H_1(W_+;\Z_2)})\pmod{8}$$
and \ref{surface} follows that
$$\chi(W)\equiv\frac{\chi(\R B)-\sigma(\C B)}{4}+\beta(q_W)\pmod{8}$$
Thus, noting that $e_A=m_1\ldots m_{s-1}m_s^2$ we get that
$$\chi(B_+)\equiv\frac{m_1\ldots m_{s-1}m_s^2}{4}+
\beta(q_{W_+}|_{H_1(V_+;\Z_2)})\pmod{8}$$
The rest of the proof is similar to that of \ref{prth}.

\subsection{Curves on an ellipsoid}
In the congruences of \ref{gen} we was able to avoid the appearance
of the complex separation with the help of the separation of $\R B$ into $B_+$
and $B_-$.
For curves of odd degree this does not work since there is no such a
separation.
Although, if $\R B$ is connected then the complex separation does not provide
the additional information and still can be avoided.
Surfaces $B_1$ and $B_2$ of the complex separation are determined by
the condition that $B_1\cup B_2=\R B$ and $\partial B_1=\partial B_2=\R A$.
The simplest case is the case of curves of odd degree on an ellipsoid.

It is well-known that the complex quadric is isomorphic to $\C P^1\times\C
P^1$,
an algebraic curve in quadric is determined by a bihomogeneous polynomial of
bidegree $(d,r)$.
If the curve is real and the quadric is an ellipsoid then $d=r$,
otherwise the curve can not be invariant under the complex conjugation
of an ellipsoid since the complex conjugation of an ellipsoid acts on
$H_2(\C P^1\times\C P^1=\Z\times\Z$ in the following way :
$conj_*(a,b)=(-b,-a)$ as it is easy to see considering the behavior of
$conj$ on the generating lines of $\C P^1\times\C P^1$.
Thus a real curve on an ellipsoid is the intersection of the ellipsoid and
a surface of degree $d$, this can be regarded as a definition of real curves
on an ellipsoid.

\subsubsection{Theorem}
\label{ell}
Let $A$ be a nonsingular real curve of bidegree $(d,d)$ on ellipsoid $B$.
Suppose that $d$ is odd.
\begin{itemize}
\begin{description}
\item[a)] If $A$ is an M-curve then
$$\chi(B_1)\equiv\chi(B_2)\equiv\frac{d^2+1}{2}\pmod{8}$$
\item[b)] If $A$ is an (M-1)-curve then
$$\chi(B_1)\equiv\frac{d^2+1}{2}\pm 1\pmod{8}$$
$$\chi(B_2)\equiv\frac{d^2+1}{2}\mp 1\pmod{8}$$
\item[c)] If $A$ is an (M-2)-curve and
$$\chi(B_1)\equiv\frac{d^2+1}{2}+4\pmod{8}$$
then $A$ is of type I
\item[d)] If $A$ is of type I then
$$\chi(B_1)\equiv\chi(B_2)\equiv 1\pmod{4}$$
\end{description}
\end{itemize}
{\em\underline{Proof}} $Dw_2(\C B)+[\R B]+[\C A]=0$ since an ellipsoid is
of type I$rel$.
Theorem 1 follows \ref{ell} now since
$$\frac{e_A}{4}+\frac{\chi(\R B)-\sigma(\C B)}{4}=\frac{d^2+1}{2}$$
and $\beta_j=0$ because of the triviality of $H_1(\R B;\Z_2)$.

\subsubsection{Low-degree curves on an ellipsoid}
Consider the application of \ref{ell} to the low-degree curves on an
ellipsoid.
Gudkov and Shustin \cite{GSh} classified real schemes
of curves of bidegree not greater then (4,4) on an ellipsoid.
To prove restrictions for such a classification it was enough to apply
the Harnack inequality and Bezout theorem.
Using the Bezout theorem avoided the extension of such a classification to
the flexible curves.

For the curves of bidegree (4,4) one can avoid the using of the Bezout theorem
using instead the analogue of the strengthened Arnold inequalities for curve
on an ellipsoid.
for the classification of the flexible curves of bidegree (3,3) the old
restrictions does not give the complete system of restrictions,
they do not restrict schemes $2\sqcup 1$$<$$2$$>$, $3\sqcup 1$$<$$1$$>$ and
$2\sqcup 1$$<$$1$$>$ (see \cite{V} for the notations).
Theorem \ref{ell} restricts these schemes and thus completes the classification
of the real schemes of flexible curves of bidegree (3,3) on an ellipsoid.

\subsubsection{Theorem}
The real scheme of a nonsingular flexible curve of bidegree (3,3) on
an ellipsoid is $1$$<1$$<$$1$$>$$>$ or $\alpha\sqcup 1$$<$$\beta$$>,
\alpha>\beta,\alpha+\beta\leq 4$.
Each of this schemes is the real scheme of some flexible (and even algebraic)
curve of bidegree (3,3) on an ellipsoid

\subsubsection{Curves of bidegree (5,5) on an ellipsoid}
Theorem \ref{ell} gives an essential restriction for the M-curves of
bidegree (5,5) on an ellipsoid, \ref{ell} leaves unrestricted 18 possible
schemes of flexible M-curves of bidegree (5,5) and 15 of them can be
realized as algebraic M-curves given by birational transformations of
the appropriate affine curves of degree 6 and another one by Viro
technique of small perturbation from the product of five plane sections
intersecting in two different points.
Thus there are only two schemes of M-curves of bidegree (5,5)
unrestricted and unconstructed, namely, $1\sqcup 1$$<$$6$$>$$\sqcup 1$$<$$8$$>$
and $1\sqcup 1$$<$$5$$>$$\sqcup 1$$<$$9$$>$.

\subsubsection{The case of bidegree (4,4)}
The classification of the real schemes of curves of bidegree (4,4) \cite{GSh}
shows that there is no congruence like \ref{ell} for curves of even degree:
the Euler characteristic of surfaces $B_1$ and $B_2$ for curves of bidegree
(4,4) can be any even number between -8 and 10.
Thus the condition that $Dw_2(\C B)+[\R B]+[\C A]=0$ is essential.

\subsubsection{Addendum, the Fiedler congruence for curves on an ellipsoid}
Let $A$ be an M-curve of bidegree $(d,d)$ on ellipsoid $B$.
Suppose that $d$ is even, the Euler characteristic of each component
of $B_1$ is even and $\chi(B_1)\equiv 2\pmod{4}$
then
$$\chi(B_2)\equiv d^2\pmod{16}$$
$$\chi(B_1)\equiv 2-d^2\pmod{16}$$
{\em\underline{Proof}} The formula of complex orientations for curves on
an ellipsoid (see \cite{Z1}) follows that the surface $V$ equal to $B_2\cup
A_+$
is $\Z_2$-homologous to zero in $\C B$ where $A_+$ is one of the components
of $\C A-\R A$.
Thus $V$ is a characteristic surface in $\C B$.
Further arguments are similar to that of \cite{F}

\subsection{Curves on cubics}
In this subsection we consider the application of Theorem 1 to
the curves of degree 2 on cubics, surfaces in $P^3$ given by cubic polynomial
(cf.\cite{Mi}) of type I$rel$.
In notations of \ref{gen} we have $q=3$, $s=2$, $m_1=3$, $m_2=2$.
Rokhlin congruence for curves on surfaces gives the complete system
of restrictions for curves of degree 2 on $M$-cubic but for restrictions
for curves of degree 2 on another cubic of type I$rel$ -- cubic diffeomorphic
to the disjoint sum of $\R P^2$ and $S^2$ we need some new tools.
Theorem 1 suffices for this purpose.
To apply Theorem 1 note that the complex separation of the disjoint cubic
consists of two surfaces -- $\R P^2$ and $S^2$.
The following is the classification of the real schemes of flexible curves
of degree 2 on a hyperboloid obtained in \cite{Mi}, for more details and
for the classifications on other real cubics see \cite{Mi}.

\subsubsection{Theorem}
Each flexible M-curves of degree 2 on non-connected cubic has one of
the following real schemes.
\begin{itemize}
\begin{description}
\item[a)] $($$<$$3\sqcup 1$$<$$1$$>$$>$$)_{\R P^2}\sqcup
($$<$$\emptyset$$>$$)_{S^2}$
\item[b)] $($$<$$1$$<$$4$$>$$>$$)_{\R P^2}\sqcup
($$<$$\emptyset$$>$$)_{S^2}$
\item[c)] $($$<$$\alpha$$>$$)_{\R P^2}\sqcup
($$<$$5-\alpha$$>$$)_{S^2}, 0\leq\alpha\leq 5$
\end{description}
\end{itemize}
Each of these 8 schemes is the real scheme of some flexible curve of degree 2
on an ellipsoid.

\end{document}